\documentclass[12pt]{amsart}

\usepackage{color}
\usepackage[centertags]{amsmath}
\usepackage[margin=1in]{geometry}
\usepackage{amsfonts}
\usepackage{amsthm}
\usepackage{newlfont}
\usepackage{amscd}
\usepackage{amsgen}
\usepackage{amssymb}
\usepackage{stmaryrd}
\usepackage{amssymb,amsmath}
\usepackage{amscd}
\usepackage{enumerate}
\usepackage{verbatim}
\usepackage{youngtab}
\newlength{\defbaselineskip} 
\theoremstyle{plain}
\newtheorem{thm}{Theorem}[section]
\newtheorem{cor}[thm]{Corollary}
\newtheorem{claim}[thm]{Claim}
\newtheorem{con}[thm]{Conjecture}
\newtheorem{df}[thm]{Definition}

\newtheorem{lema}[thm]{Lemma}
\newtheorem{obs}[thm]{Proposition}
\newtheorem{exm}[thm]{Example}
\newtheorem{fact}[thm]{Fact}

\newtheorem{rem}[thm]{Remark}
\newtheorem{pr}{Algorithm}
\theoremstyle{definition} 
\theoremstyle{definition} \newtheorem{ex}{Example}[section] %

 \numberwithin{equation}{section}

\DeclareMathOperator{\sgn}{sgn}

\def\ob{\begin{obs}}
\def\kob{\end{obs}}
\def\fa{\begin{fact}}
\def\kfa{\end{fact}}
\def\dow{\begin{proof}}
\def\kdow{\end{proof}}

\def\tw{\begin{thm}}
\def\ktw{\end{thm}}
\def\hip{\begin{con}}
\def\khip{\end{con}}
\def\lem{\begin{lema}}
\def\klem{\end{lema}}
\def\cl{\begin{claim}}
\def\kcl{\end{claim}}
\def\ex{\begin{exm}}
\def\prog{\begin{pr}}
\def\kprog{\end{pr}}
\def\wn{\begin{cor}}
\def\kwn{\end{cor}}
\def\uwa{\begin{rem}}
\def\kuwa{\end{rem}}
\def\kex{\end{exm}}
\def\dfi{\begin{df}}
\def\kdfi{\end{df}}

\def\CC{\mathbf{C}}

\begin{document}
\title{Effective constructions in plethysms  and
Weintraub's conjecture}
%\author{Laurent Manivel, Mateusz Micha\l ek}
\author[L. Manivel]{Laurent Manivel}
\address{Laurent Manivel\newline
Institut Fourier,
UMR 5584 Universit\'e de Grenoble/CNRS,
BP 74, 38402 Saint-Martin d'H\`eres, France}
\email{{\tt Laurent.Manivel@ujf-grenoble.fr}}
\author[M. Micha\l ek]{Mateusz Micha\l ek}
\address{Mateusz Micha\l ek\newline
Max Planck Institute for Mathematics,
Vivatsgasse 7,
53111 Bonn,
Germany
}
\email{{\tt wajcha2@poczta.onet.pl}}
\thanks{The second author is supported by the NCN grant UMO-2011/01/N/ST1/05424}
\maketitle

%\dfi[Young diagram $(a_1,\dots,a_l)$]
%Assume that $a_1\geq\dots\geq a_l$ are positive integers. The Young diagram $(a_1,\dots,a_l)$ i
%s a diagram with $l$ columns. The length of the $i$-th column is $a_i$.
%\kdfi
\begin{abstract}
 We give a short proof of Weintraub's conjecture \cite{we},
first proved in \cite{bci}, by constructing explicit highest
weight vectors in the plethysms $S^p(\wedge^{2q}W)$.
\end{abstract}

\section{Introduction}

Plethysm is one of the most basic operations on symmetric functions.
It was introduced by Littlewood along his fundamental work on group
representations, and it has remained notably difficult to understand and
to compute (see e.g. \cite[I.8]{mcd}, and also \cite{lor} for a recent
overview and more references).
In the language of representation theory, plethysm is defined as the
composition of Schur functors. Even the very special cases of
compositions of symmetric or skew-symmetric powers seem
completely out of reach of our current tools.

Among the few known general properties of plethysm, one was first
observed in low degrees and then conjectured by Weintraub \cite{we}.
The claim is that for any partition $\lambda = (\lambda_1\ge \ldots\ge
\lambda_p\ge 0)$, with only {\it even} parts, whose sum is $2pq$,
then the Schur module $S_{\lambda}W$ appears in the composition of
symmetric powers $S^p(S^{2q}W)$ with non zero multiplicity. (Here
the complex vector space $W$ is supposed to be of dimension at least
equal to the number of non zero parts in $\lambda$, and then this
multiplicity does not depend on it. We refer to \cite{fh} for the
definition of Schur modules and basic facts on the representation
theory of $GL(W)$.) An asymptotic version  was
established in \cite{man}, but the conjecture in its full generality
was first proved in \cite{bci}, using ideas and methods from quantum
information theory.

The main purpose of our paper is to give a different, more traditional
proof of Weintraub's conjecture. Another motivation of our work being
related to the Pl\"ucker embeddings of Grassmannians, we rather
consider symmetric powers of wedge powers. This doesn't make any
difference regarding Weintraub's conjecture because of the ``duality''
$$S^p(S^{2q}W) = S^p(\wedge^{2q}W)^*$$
(see \cite[I.8, Ex. 1]{mcd} for the corresponding statement in terms
of symmetric functions).
This duality statement should be understood as follows: for any partition
$\lambda$ such that  $S_{\lambda}W$ appears inside $S^p(S^{2q}W)$, then
$S_{\lambda^*}W$ appears inside  $S^p(\wedge^{2q}W)$ with the same
multiplicity, where the dual partition $\lambda^*$ is defined by
$\lambda^*_i=\# \{j, \lambda_j\ge i\}$. (Once again, in order for
that statement to be correct one needs to suppose that the dimension
of $W$ is large enough, namely larger or equal to $2pq$.)

The proof we give of Weintraub's conjecture consists in providing an
explicit construction of a highest weight vector of weight $\lambda^*$
inside $S^p(\wedge^{2q}W)$. We first give in section 3.1 an algorithm
to construct a special weight vector inside $(\wedge^{2q}W)^{\otimes p}$.
In 3.2 we check that it is a highest weight vector, hence that it
contributes to the multiplicity of $S_{\lambda^*}W$ in this tensor product.
Finally in 3.3 we prove that its symmetrization is non zero inside
$S^p(\wedge^{2q}W)$, which implies the conjecture. Section 2 is essentially
a warm-up. In 2.1 we explain how to construct basis of highest weight vectors
inside $(\wedge^{2q}W)^{\otimes p}$ using similar, but even simpler
ideas. In 2.2 we take a slightly different perspective on these
highest weight vectors and deduce some consequences on asymptotic
multiplicities.

\section{Highest weight vectors in tensor products}

\subsection{Highest weight vectors in tensor powers of wedge products}
Let $a_{k,d}(\lambda)$ denote the multiplicity of the irreducible $GL(W)$-module
$S_{\lambda^*} W$ inside $(\bigwedge^k W)^{\otimes d}$. This multiplicity does not depend
on the dimension of $W$, provided that $\dim W \ge \ell(\lambda)$, the number of non-zero
parts of $\lambda$, which we will
always suppose. We can calculate $a_{k,d}(\lambda)$ using Pieri's rule: it is
equal to the number of tableaux $T$ of shape $\lambda^*$ and weight $(k^d)$,
which are increasing on rows and non decreasing on columns. (We refer to
\cite{mcd} for Pieri's rule and the language of tableaux, see also \cite{fh}.)

To each such tableau $T$ we will associate a highest weight vector $w_T$ of
weight $\lambda^*$ in $(\bigwedge^k W)^{\otimes d}$. We will show that
the vectors $w_T$, when $T$ varies, form a basis of the highest weight space.

Let $h_a^b$ be the number of entries  equal to $b$ in the $a$-th column of $T$.
Let $f_a^b=\sum_{i=1}^b h_a^i$ be the number of entries in the $a$-th column that
are less or equal to $b$. Note that $h_a^b=0$ for $a>b$ and $f_1^1=h_1^1=k$.
Each of these sequences completely encodes the tableau  $T$.
%Let $\lambda$ be a Young diagram with columns of lengths $\lambda_1\geq \lambda_2\geq \dots$.
%Let $h$ be the sequence of all $h_a^b$ that encodes a given decomposition according
%to Pieri's rule. We order the sequence $h$ lexicographically.
%\ex
%Continuing the previous example we have $h_1^1=3$, $h_1^2=1$ and $h_2^2=2$.
%\kex
%To each column of the Young diagram of $\lambda$ of length $\lambda_i$ we associate
%a group of permutations $S_{\lambda_i}$. We are ready to state the definition
%of vectors associated to a given decomposition $h$ of the diagram
%$\lambda=(\lambda_1,\dots,\lambda_l)$. Let $e_i$ be the fixed basis of $W$.
\dfi[Vectors $w_T$, simple tensors $t_{\gamma_1,\dots,\gamma_l}$]\label{baza}
%Let us fix a Young diagram filled according to Pieri's rule, represented by $h$.
Let $e_1,\ldots ,e_N$ be the fixed basis of $W$.
For each collection of permutations $\gamma_i\in S_{\lambda_i}$, let
$$t_{\gamma_1,\dots,\gamma_l}:=\bigotimes_{g=1}^d(\bigwedge_{s=1}^{l} e_{\gamma_s(f^{g-1}_s+1)}
\wedge\dots\wedge e_{\gamma_s(f^g_s)})\;\in (\bigwedge^k W)^{\otimes d}$$
Then we associate to $T$ the vector $w_T$ given by
$$w_T:=\frac{1}{\prod h_i^j!}\sum_{\gamma_i\in S_{\lambda_i}}\large(\prod_i\sgn(\gamma_i)\large)
t_{\gamma_1,\dots,\gamma_l}.$$
\kdfi
We divide by the normalizing factor $\prod h_i^j!$ because if two permutations differ
only on entries that appear in the same tensor product, then these permutations define the
same simple tensor. In particular $w_T$ has only integer coefficients.
\ex
The multiplicity of $S_{2,2,1,1}W$ inside $(\bigwedge^3 W)^{\otimes 2}$ equals $1$.
The unique suitable tableau $T$ is
$${\tiny\young(12,12,1,2)}.$$
We have $h_1^1=3$, $h_1^2=1$ and $h_2^2=2$. The corresponding vector is
\begin{eqnarray*}
w_T & = & \frac{1}{3!\cdot 2!}\sum_{\gamma_1\in S_4,\gamma_2\in S_2}\sgn(\gamma_1)\sgn(\gamma_2)
(e_{\gamma_1(1)}\wedge e_{\gamma_1(2)}\wedge e_{\gamma_1(3)})\otimes(e_{\gamma_1(4)}
\wedge e_{\gamma_2(1)}\wedge e_{\gamma_2(2)}) \\
 & = & (e_1\wedge e_2\wedge e_3)\otimes(e_4\wedge e_1\wedge e_2)-(e_1\wedge e_2\wedge e_4)
\otimes(e_3\wedge e_1\wedge e_2).
\end{eqnarray*}
\kex
%\ex
%Consider $(W\wedge W)\otimes (W\wedge W)$. According to Pieri's rule this representation contains a subrepresentation corresponding to Young diagram $(3,1)$. This is given by $h^2_2=h^2_1=1$. The corresponding vector $w_h$ is equal to:
%$$(e_1\wedge e_2)\otimes (e_3\wedge e_1)-(e_1\wedge e_3)\otimes (e_2\wedge e_1).$$
%\kex
\begin{comment}
The following definition will be useful in forthcoming sections, especially while considering limit values of multiplicities.
\dfi[Type $A$, type $B$]
For a given $w_h$ and $1\leq i\leq d$ we say that the $i$-th wedge product is of type $A$ if it depends on $\sigma_1$ only.
For example according to Pieri's rule the first wedge product is always of type $A$ for any $w_h$. Products not of type $A$ are called of type $B$.
\kdfi
Products of type $A$ correspond in Pieri's rule to blocks entirely contained in the first column.
\end{comment}

\ob\label{basis}
The vectors $w_T$, for $T$ of shape $\lambda^*$, form a basis of the highest weight space
of weight $\lambda^*$ inside $(\bigwedge^k W)^{\otimes d}$.
\kob

\proof The fact that each $w_T$ is a highest weight vector is checked by a routine
computation. There remains to show that they are linearly independent. We need a definition.
\dfi[Vectors $r_h$]\label{bazadualna}
Let  $(e_1^*,\ldots ,e_N^*)$ be the dual basis to $(e_1,\ldots ,e_N)$.
To a tableau $T$ encoded by the sequences $h$ or $f$ we associate
$$r_T=\bigotimes_{g=1}^d (\bigwedge_{s=1}^{l(\lambda)} e_{f^{g-1}_s+1}^*\wedge\dots\wedge
e_{f^g_s}^*)\;\in (\bigwedge^k W^*)^{\otimes d}.$$
\kdfi
We order the tableaux as follows. If $T$ and $T'$ are two tableaux, encoded
by two sequences $h$ and $h'$, we let $T'>T$ if and only if $h'>_{lex}h$, where
for the lexicographic order the sequences $h$ and $h'$ are read column after
column. The following easy lemma is left to the reader.
\lem\label{dual}
For $T'>T$, $r_{T'}(w_T)=0$. Moreover $r_{T}(w_T)=1$.
\klem

This clearly implies that the vectors $w_T$ are independent. \qed
%For $h>h'$ the value of $r_{h'}(w_{h})$ does not have to be zero, although sometimes it is.
%In such cases the computations are particulary easy.

\medskip\noindent {\it Remark}. Of course the same method would allow to produce
basis of highest weight vectors in any tensor product of wedge powers, and can be also
adapted to symmetric powers.

\subsection{Highest weight vectors and asymptotic multiplicities}
In this section we consider this question under the slightly different
perspective of computing asymptotic multiplicities: that is, we consider
partitions with a varying first row, the remaining part being fixed. It
has been observed that the corresponding multiplicities inside plethysms
of symmetric powers for example, are non-decreasing functions of the
exponents, and becomes eventually constant \cite{man}. We will give
a simple interpretation of certain of the asymptotic multiplicities, those
multiplicities that are obtained when the exponents are large enough.
We focus on the case of symmetric powers, which is slightly simpler.

\medskip
By Pieri's rule, the multiplicity $s_{k,d}(\lambda)$
of $S_{kd-|\lambda|, \lambda}V$ inside $(S^kV)^{\otimes d}$
is equal to the number of semistandard tableaux $T$ of shape $(kd-|\lambda|, \lambda)$,
with  $k$ entries equal to $i$ for each $1\le i\le d$. Such a tableau is completely
determined by the part below the first row, whose entries are all bigger than one.
Substracting one to each entry we get a semistandard tableau $S$ of shape $\lambda$,
with entries between $1$ and $d-1$. Conversely, if $S$ is such a tableau, and if the
number of occurrences $e_i(S)$ of each entry $i$ in $S$ is not greater than $k$, we can
recover $T$ by adding one to each entry of $S$, and then $k$ one's above the first row,
$k-e_1(S)$ two's, etc... The resulting tableau $T$ is certainly semistandard if
$\lambda_1\le k$. Under this hypothesis we therefore get a bijection between
two types of tableaux. Observe that the number of tableaux $S$ is equal to the
dimension of the Schur power $S_\lambda \CC^{d-1}$, hence the equality
$$s_{k,d}(\lambda)=\dim S_\lambda \CC^{d-1}\qquad \mathrm{for}\quad k\ge \lambda_1.$$
We will give a more precise version of this equality, as an identity between representations
of the symmetric group. Recall that the fundamental representation of the symmetric
group $S_d$, denoted by $[d-1,1]$, is obtained by permutating coordinates $x_1,\ldots ,x_d$
in the hyperplane of $\CC^d$ of equation $x_1+\cdots +x_d=0$.

\begin{obs}
For $k\ge \lambda_1$ and $\dim V>\lambda_1^*$, there is an isomorphism
of $S_d$-modules
$$Hom_{GL(V)}(S_{kd-|\lambda|, \lambda}V,(S^kV)^{\otimes d})\simeq S_\lambda[d-1,1].$$
\end{obs}

As a consequence of Schur-Weyl duality we have the identity
$$(S^kV)^{\otimes d}=\bigoplus_{|\mu|=d}S_\mu(S^kV)\otimes [\mu]$$
of $GL(V)\times S_d$-modules. Hence the following corollary,
which can also be extracted from \cite[Corollary 5.3]{br}:

\begin{cor}
Let $k\ge \lambda_1$ and $\dim V>\lambda_1^*$. Then for any partition $\mu$
of size $d$, the multiplicity of
$S_{kd-|\lambda|, \lambda}V$ inside $S_\mu(S^kV)$
is equal to the multiplicity of $[\mu]$ inside $S_\lambda[d-1,1]$.
\end{cor}

\noindent {\it Proof of the Proposition}. We choose a basis $e_0,\ldots ,e_N$ of $V$.
Then the  space of $GL(V)$-equivariant morphisms $Hom_{GL(V)}(S_{kd-|\lambda|,
\lambda}V,(S^kV)^{\otimes d})$ can be identified with the space of highest
weight vectors of weight $(kd-|\lambda|, \lambda)$ inside $(S^kV)^{\otimes d}$.
A basis of the latter space is given by monomials $e_0^{k-|\alpha_1|}e^{\alpha_1}
\otimes\cdots\otimes e_0^{k-|\alpha_d|}e^{\alpha_d}$, where each $\alpha_i$
is a sequence of $N$ integers, with sum $|\alpha_i|\le k$. The subspace
$(S^kV)^{\otimes d}_{(kd-|\lambda|, \lambda)}$
of vectors of weight $(kd-|\lambda|, \lambda)$ is isomorphic with $S^{\lambda_1}\CC^d
\otimes\cdots\otimes S^{\lambda_N}\CC^d$. If $f_1,\ldots ,f_d$ is a basis of
$\CC^d$, an explicit isomorphism $\theta$ is obtained by sending the monomial
$e_0^{k-|\alpha_1|}e^{\alpha_1}
\otimes\cdots\otimes e_0^{k-|\alpha_d|}e^{\alpha_d}$ to the monomial
$$f_1^{\alpha_{1,1}}\cdots f_d^{\alpha_{d,1}}\otimes \cdots\otimes
f_1^{\alpha_{1,N}}\cdots f_d^{\alpha_{d,N}}.$$
This identification is compatible with the action of the symmetric group,
if $S_d$ acts on $\CC^d$ by permuting the $f_j$'s.

Now, vectors in $(S^kV)^{\otimes d}_{(kd-|\lambda|, \lambda)}$ are highest
weight vectors if and only if they are killed by each of the endomorphisms
induced  on $(S^kV)^{\otimes d}$ by the endomorphisms $X_i$ of $V$,
with $0\le i\le N-1$, that sends $e_{i+1}$ to $e_i$ and any other $e_j$ to zero.
%Of course one can restrict to $i\le \ell(\lambda)$.

Under the identification given by the isomorphism $\theta$,
the action of $X_i$ for $1\le i\le d-1$ is easily seen to
coincide with the natural morphism
$$Y_i : S^{\lambda_1}\CC^d\otimes\cdots\otimes S^{\lambda_i}\CC^d\otimes S^{\lambda_{i+1}}\CC^d
\otimes\cdots\rightarrow
S^{\lambda_1}\CC^d\otimes\cdots\otimes S^{\lambda_i+1}\CC^d\otimes S^{\lambda_{i+1}-1}\CC^d
\otimes\cdots $$
Moreover the action of $X_0$ is given by the morphism
$$Y_0 : S^{\lambda_1}\CC^d\otimes S^{\lambda_2}\CC^d
\otimes\cdots\otimes S^{\lambda_N}\CC^d\rightarrow
S^{\lambda_1-1}\CC^d\otimes S^{\lambda_2}\CC^d
\otimes\cdots\otimes S^{\lambda_N}\CC^d,$$
where the map $S^{\lambda_1}\CC^d\rightarrow S^{\lambda_1-1}\CC^d$ is induced
by the linear form $u$ sending each $f_i$ to $1$.

We conclude that the space of highest weight vectors in
$(S^kV)^{\otimes d}_{(kd-|\lambda|, \lambda)}$
can be identified with $Ker(Y_0)\cap Ker(Y_1)\cap\cdots\cap Ker(Y_{d-1})
\subset S^{\lambda_1}\CC^d\otimes\cdots\otimes S^{\lambda_N}\CC^d$.
But the intersection $Ker(Y_1)\cap\cdots\cap Ker(Y_{d-1})=S_\lambda\CC^d$, and
$Ker(Y_0)\cap S_\lambda\CC^d = S_\lambda Ker (u)\simeq S_\lambda[d-1,1]$.

Indeed, for any hyperplane $K$ of $\CC^d$ the kernel of the map
$S^p\CC^d\otimes S^q\CC^d\rightarrow S^{p+1}\CC^d\otimes S^{q-1}\CC^d$
restricted to $S^pK\otimes S^q\CC^d$ is easily seen to be contained in
$S^pK\otimes S^qK$. By induction, we deduce that $Ker(Y_0)\cap S_\lambda\CC^d$
is contained in the intersection of the kernels of $Y_1, \ldots , Y_i$
restricted to  $S^{\lambda_1}Ker(u)\otimes \cdots\otimes S^{\lambda_i}Ker(u)\otimes
S^{\lambda_{i+1}}\CC^d\otimes\cdots\otimes S^{\lambda_N}\CC^d$.
The final case $i=N$ yields the claim. \qed

\medskip
The decomposition of $S_\lambda [d-1,1]$ is a difficult problem.
It is known \cite[Lemma 7.5]{br} that
$$\wedge^i[d-1,1]=[d-i,1^i].$$
With the previous corollary this implies that $S_\mu(S^kV)$ can contain
an irreducible component of hook shape only if $\mu$ is itself a hook.
(In fact the corollary implies this only asymptotically, but since
multiplicities are known to be non decreasing functions of both
exponents, see \cite{man}, the general statement follows.) This result
first appeared in \cite{lar}.

One way to proceed in general in order to compute $S_\lambda [d-1,1]$,
would be to express the Schur functor
$S_\lambda$ in terms of exterior powers only using the Giambelli
formula. It could then be computed
by induction if we knew how to decompose tensor products by $[d-i,1^i]$.
By Schur-Weyl duality this amounts to computing $S_{d-i,1^i}(A\otimes B)$
in terms of Schur powers of the two vector spaces $A$ and $B$ (of large
enough dimensions).
But then we can use the fact that $S_{d-i,1^i}=\oplus_{j\ge 0}(-1)^j
S^{d-i+j}\otimes \wedge^{i-j}$ to reduce to a computation involving
only Littlewood-Richardson coefficients. For example, if we define
$\langle\nu\rangle$ to be the representation $[d-|\nu|,\nu]$ when it
makes sense, and zero otherwise, we get
$$\begin{array}{rcl}
   S_{2,1^i}\langle 1\rangle
 & = & \langle 2,1^i\rangle\oplus  \langle 2, 1^{i-1}\rangle\oplus
\langle 1^{i+1}\rangle\oplus  \langle 1^i\rangle, \\
% &  & \\
S_{3,1^i}\langle 1\rangle
 & = & \langle 3,1^i\rangle\oplus  \langle 3, 1^{i-1}\rangle\oplus
\langle 2^2,1^{i-2}\rangle\oplus  2\langle 2,1^i\rangle\oplus \\
 & & \oplus
\langle 2^2,1^{i-3}\rangle\oplus  3\langle 2, 1^{i-1}\rangle\oplus
3\langle 1^{i+1}\rangle\oplus  \langle 2,1^{i-2}\rangle\oplus
 2\langle 1^i\rangle.
 \end{array}$$

%Using Young symmetrizer, the vectors $w_h$ can be projected to $S_{\mu}(\bigwedge^k)$ for
%any $\mu$ of weight $d$. In this way one obtains a system of generators of the highest
%weight space of the isotypic component corresponding to $\lambda$ in the plethysm.
%In fact for small $\mu$ the projection matrix given by the Young symmetrizer in the
%described basis can be explicitly computed. In particular, the trace of such a matrix
%equals the multiplicity
% of the corresponding component in the plethysm. Moreover the multiplicity is equal
%to zero if and only if all vectors $w_h$ project to zero. Thus they form good candidates
%to prove that the multiplicity is nonzero. This approach will be used in the following
%section. Although we do not project vectors $w_h$, the presented construction is very similar.
\medskip

\section{Weintraub's conjecture}\label{Weintraub}
In this section we explain our constructive proof of  Weintraub's conjecture \cite{we}.
\tw
Suppose that $k$ is even.
Consider an even partition $\lambda$ of weight $dk$.
Then the  multiplicity of $S_{\lambda^*}W$ in $S^d(\bigwedge^k W)$ is positive.
\ktw

We will explicitly construct a vector in $S^d(\bigwedge^k W)$ that is a highest
weight vector of weight $\lambda^*$. For this we will proceed in three steps. First
we will construct a special vector $P$ inside $(\bigwedge^kW)^{\otimes d}$.
Then we will show that $P$ is a highest weight vector of weight $\lambda^*$.
Finally, we will show that the projection of $P$ to
$S^d(\bigwedge^kW)$ is nonzero.

\subsection{Construction of a special vector}
We fix a partition $\lambda$ with only even parts
$\lambda_1\geq\dots\geq\lambda_l$. We will construct a vector $P$ of weight
$\lambda^*$ inside
$(\bigwedge^kW)^{\otimes d}$ as a combination of simple tensors. A
simple tensor can be represented by a rectangular tableau of size $k\times d$,
each column
\[\young({{a_1}},{{a_2}},\vdots,{{a_k}}).\]
representing a product of basis vectors $e_{a_1}\wedge\dots\wedge e_{a_k}$.
Note that we can freely permute the entries in a same column: this will
affect the simple tensor only by  a sign.

The vector $P$ will be constructed by an algorithm that fills the entries of
the rectangle $d\times k$, that we denote by $R(k,d)$, to obtain a tableau $T$
(or rather a combination of
tableaux indexed by permutations). Each entry filled in $T$ will correspond
bijectively to a box in (the Young diagram of) $\lambda^*$. Hence in the algorithm,
each time we fill
an entry in $T$, we also cross out the corresponding box in $\lambda^*$.

After each step we will get a partial tableau $T'$, and the part of the rectangle
$R(k,d)$ that will remain to be filled will be a subrectangle $R'\cong R(k',d')$
in the lower right corner of $R(k,d)$. Let
\begin{itemize}
% \item $d'$ (resp. $k'$) the number of columns (resp. rows) that remain to be filled;
 \item $m'=k-k'$; it will always be even;
  \item $l'$ be the number of columns of $\lambda$ with the entry in the $m'+1$-st row not crossed out.
\end{itemize}
It is very important to keep in mind that throughout the algorithm we always have $l'\leq d'$.
For $l'>0$ the number $m'$ will be equal to the number of rows of $\lambda^*$ that are completely
crossed out.

Each new step will depend on a specific column of $\lambda^*$, the leftmost
column among those that have not been completely crossed out yet.
We will denote by:
\begin{itemize}
 \item $o'$ the index of this column;
 \item $h'$ the number of boxes in that column that have already
been crossed out;
 \item $j'$ the number of boxes in that column that have not already
been crossed out.
\end{itemize}

The algorithm has three possible steps.

\medskip\noindent {\bf Step $A$}. This step applies when $l'=d'$. This means
that the number of columns that we still have to deal with in $\lambda^*$ is equal
to the number of columns in $R'$.
Then we fill the two top rows of $T'$: the first one with $m'+1$ and the second
one with $m'+2$. The corresponding entries will be called {\it frozen}.
We cross out in $\lambda^*$ the first two entries of each column
that we still have to deal with.

After this step $d'$ remains unchanged (unless we filled the whole rectangle)
while $l'$ might decrease or remains unchanged. In particular we still have
$l'\leq d'$.

\smallskip The two other possible steps apply when $l'<d'$. They will fill the
leftmost column of $R'$. In particular after each of these steps $d'$ decreases
by one while $l'$ might decrease or remain unchanged. In particular the relation
$l'\leq d'$ is preserved.

\medskip\noindent {\bf Step $B$}. This step applies when $l'<d'$ and $j'\ge k'$.
Then we fill the leftmost column of $R'$ with $\sigma_{o'}(h'+1),\sigma_{o'}(h'+2),\dots,
\sigma_{o'}(h'+k')$ starting from the top. In $\lambda^*$ we cross-out the $k'$
topmost boxes.

\medskip\noindent {\bf Step $C$}. This step applies when $l'<d'$ and $j'<k'$.
This means that we have less entries left in the column $o'$  than entries
to fill in a column of $R'$, so we will need to pass to another column of $\lambda^*$.
First we deal with the column $o'$. In the leftmost column of $R'$
we fill the boxes with $\sigma_{o'}(h'+1),\dots,
\sigma_{o'}(h'+j')$ starting from the top. Correspondingly, we completely
cross out the column $o'$. Then we pass to the next column of $\lambda^*$
where we will cross out the missing number of boxes.
Note that this column has exactly $m'$ boxes already crossed out. We complete the
leftmost column of $R'$ by $\sigma_{o'+1}(m'+1),\dots,
\sigma_{o'+1}(m'+k'-j')$. By Claim \ref{niemadwoch} we can do that without
having to go to the next column of $\lambda^*$.
\medskip

\dfi
The vector $P$ is the sum of all the simple tensors associated to the
tableaux produced by the algorithm, weighted by the product of the
signs of the permutations involved. Note that $P$ is certainly a weight
vector of weight $\lambda^*$.
\kdfi

We insist on the fact that the vector $P$ is represented by a single
tableau, constructed by the algorithm. In this tableau certain entries
are frozen. All the other entries are affected by a permutation that
depends only on the column of the corresponding box in $\lambda$.
\medskip

\ex
Consider $\lambda^*=(4,4,3,3,3,3)$, the dual of $\lambda=(6,6,6,2)$,
%\[\yng(5,5,4,4,2,2,1,1)\]
 and the tensor product $(\bigwedge^4W)^{\otimes}$. Let us apply the algorithm.

\medskip
At the beginning we have $l'=4<d'=d=5$. Thus we do not perform step A.
Since $\lambda_1=6\geq k'=k=4$ we apply step B and we obtain:
 \vskip 8mm
$$\hskip -8.2cm T=$$\vskip -1.4cm { $$\hskip -3.7cm{\fontsize 11\young({{\sigma(1)}}{{}}{{}}{{}}{{}},{{\sigma(2)}}{{}}{{}}{{}}{{}},{{\sigma(3)}}{{}}{{}}{{}}{{}},{{\sigma(4)}}{{}}{{}}{{}}{{}})}$$}
 \vskip -2.8cm$$\hskip 2cm \lambda^* =$$\vskip -1.8cm{$$\hskip 5.6cm\young({{\times}}{{}}{{}}{{}},{{\times}}{{}}{{}}{{}},{{\times}}{{}}{{}},{{\times}}{{}}{{}},{{}}{{}}{{}},{{}}{{}}{{}}).$$}

We get $l'=3< d'=4$. As there are only two entries left in the first column of $\lambda$ we have to
apply step C. We get:
\vskip 5mm$$\hskip -8.2cm T=$$\vskip -1.4cm { $$\hskip -3.7cm{\fontsize 11\young({{\sigma(1)}}{{\sigma(5)}}{{}}{{}}{{}},{{\sigma(2)}}{{\sigma(6)}}{{}}{{}}{{}},{{\sigma(3)}}{{\tau(1)}}{{}}{{}}{{}},{{\sigma(4)}}{{\tau(2)}}{{}}{{}}{{}})}$$}
 \vskip -2.8cm$$\hskip 2cm \lambda^*=$$\vskip -1.8cm{$$\hskip 5.6cm\young({{\times}}{{\times}}{{}}{{}},{{\times}}{{\times}}{{}}{{}},{{\times}}{{}}{{}},{{\times}}{{}}{{}},{{\times}}{{}}{{}},{{\times}}{{}}{{}}).$$}

Now $l'=2<d'=3$ and there are four entries left in the second column of $\lambda$.
So we apply step B and obtain:
\vskip 5mm$$\hskip -8.2cm T=$$\vskip -1.4cm{ $$\hskip -3.7cm{\fontsize 11\young({{\sigma(1)}}{{\sigma(5)}}{{\tau(3)}}{{}}{{}},{{\sigma(2)}}{{\sigma(6)}}{{\tau(4)}}{{}}{{}},{{\sigma(3)}}{{\tau(1)}}{{\tau(5)}}{{}}{{}},{{\sigma(4)}}{{\tau(2)}}{{\tau(6)}}{{}}{{}})}$$}
 \vskip -3.2cm$$\hskip 2cm \lambda^*=$$\vskip -1.8cm{$$\hskip 5.6cm\young({{\times}}{{\times}}{{}}{{}},{{\times}}{{\times}}{{}}{{}},{{\times}}{{\times}}{{}},{{\times}}{{\times}}{{}},{{\times}}{{\times}}{{}},{{\times}}{{\times}}{{}}).$$}

Now we get $l'=d'=2$, so we apply step A to obtain:
\vskip 5mm$$\hskip -8.2cm T=$$\vskip -1.4cm{ $$\hskip -3.7cm{\fontsize 11\young({{\sigma(1)}}{{\sigma(5)}}{{\tau(3)}}{{1}}{{1}},{{\sigma(2)}}{{\sigma(6)}}{{\tau(4)}}{{2}}{{2}},{{\sigma(3)}}{{\tau(1)}}{{\tau(5)}}{{}}{{}},{{\sigma(4)}}{{\tau(2)}}{{\tau(6)}}{{}}{{}})}$$}
 \vskip -3.2cm$$\hskip 2cm \lambda^*=$$\vskip -1.8cm{$$\hskip 5.6cm\young({{\times}}{{\times}}{{\times}}{{\times}},{{\times}}{{\times}}{{\times}}{{\times}},{{\times}}{{\times}}{{}},{{\times}}{{\times}}{{}},{{\times}}{{\times}}{{}},{{\times}}{{\times}}{{}}).$$}

%. We will apply twice step B. Consider a permutation $\delta$ of the set $\{3,4,5,6\}$. We have:
%\vskip 5mm$$\hskip -8.2cm T=$$\vskip -1.4cm { $$\hskip -3.7cm{\fontsize 11\young({{\sigma(1)}}{{\sigma(5)}}{{1}}{{1}}{{1}}{{1}},{{\sigma(2)}}{{\sigma(6)}}{{2}}{{2}}{{2}}{{2}},{{\sigma(3)}}{{\sigma(7)}}{{\delta(3)}}{{\delta(5)}}{{}}{{}},{{\sigma(4)}}{{\sigma(8)}}{{\delta(4)}}{{\delta(6)}}{{}}{{}})}$$}
% \vskip -2.8cm$$\hskip 2cm \lambda=$$\vskip -1.8cm{$$\hskip 5.6cm\young({{\times}}{{\times}}{{\times}}{{\times}}{{\times}},{{\times}}{{\times}}{{\times}}{{\times}}{{\times}},{{\times}}{{\times}}{{}}{{}},{{\times}}{{\times}}{{}}{{}},{{\times}}{{\times}},{{\times}}{{\times}},{{\times}},{{\times}}).$$}

Now $l'=1<d'=2$, so we apply step B. Finally, as $l'=0<d'=1$ we apply step B and get:
\vskip 1mm$$\hskip -4.2cm T=$$\vskip -1.4cm { $${\fontsize 11\young({{\sigma(1)}}{{\sigma(5)}}{{\tau(3)}}{{1}}{{1}},{{\sigma(2)}}{{\sigma(6)}}{{\tau(4)}}{{2}}{{2}},{{\sigma(3)}}{{\tau(1)}}{{\tau(5)}}{{\delta(3)}}{{\delta(5)}},{{\sigma(4)}}{{\tau(2)}}{{\tau(6)}}{{\delta(4)}}{{\delta(6)}})}$$}
%\vskip 5mm$$\hskip -8.2cm T=$$\vskip -1.4cm { $$\hskip -3.7cm{\fontsize 11\young({{\sigma(1)}}{{\sigma(5)}}{{1}}{{1}}{{1}}{{1}},{{\sigma(2)}}{{\sigma(6)}}{{2}}{{2}}{{2}}{{2}},{{\sigma(3)}}{{\sigma(7)}}{{\delta(3)}}{{\delta(5)}}{{3}}{{3}},{{\sigma(4)}}{{\sigma(8)}}{{\delta(4)}}{{\delta(6)}}{{4}}{{4}}).}$$}
\kex
\medskip

\subsection{The vector $P$ is a highest weight vector}
We need to show  that $P$ is killed by each of the operator $X_j$ that sends
$e_{j+1}$ to $e_j$ and any other $e_i$ to zero.
Let us consider the possible occurrences of $e_{j+1}$ in $P$
and how they came to appear when we applied the algorithm.

\smallskip
If $e_{j+1}$ has been produced by step A, than it is frozen on the
corresponding column of $\lambda$ and $j$ is also frozen on this column.
So applying $X_j$ we certainly get zero.

\smallskip
If $e_{j+1}$ has been produced by step B or C, then $j$ could be already
frozen in the corresponding column of $\lambda$ and the same argument applies.
Otherwise, $j$ and $j+1$ are affected by the same permutation. The terms
for which $e_j$ and $e_{j+1}$ appear in the same column of the tableau are certainly killed
by $X_j$; while the terms for which  $e_j$ and $e_{j+1}$
appear in different columns come in pairs, just be switching $e_j$ and $e_{j+1}$;
these terms have different signs but the same image by $X_j$, so their
contributions cancel out.

\subsection{The symmetric projection of $P$ is non zero}
We have not yet proved that $P$ is non zero. We will show directly that
its projection to $S^d(\bigwedge^kW)$ is non zero.

\dfi[Vector $Q$]
Consider the simple tensor in $P$ that is obtained by taking all the permutations  equal
to the identity. By symmetrizing this simple tensor we get the vector $Q$ in $S^d(\bigwedge^kW)$.
\kdfi

We will proceed as follows. We will first show that the vector $Q$ is nonzero. Then we will
prove that each time $Q$ is obtained as the symmetrization of a simple tensor in $P$,
it comes with a positive sign. This will certainly imply the claim.

\medskip\noindent {\bf $Q$ is non zero}.
What we need to check that if we take all permutations in the
definition of $P$ equal to the identity, there is no repetition
in any column of the resulting rectangular tableau.

No repetition can come from the frozen variables.
Let us observe that in step $C$ we decrease $l'$, thus after this step the
strict inequality $l'<d'$ holds. It follows that step $A$ cannot follow
immediately after $C$.
Note that after step $B$ we have $l'\geq m'+k'=k$. Thus, each time we apply
step $A$ we have $l'\geq k$. Of course the frozen entries are always less or equal
to $k$. Hence the frozen entries cannot coincide with entries filled in step $B$ or $C$.

Moreover no
repetition can appear when we apply step B since in this case
only one column of $\lambda$ is involved. There remains to
consider step C in more detail. It seems to be the right moment
to prove

\cl\label{niemadwoch}
Let $y_1=\lambda_{o'}-h'$ be the number of uncrossed boxes in the column
that we are dealing with, and $y_2=\lambda_{o'+1}-m'$ the number of uncrossed
boxes in the next column. Then $y_1+y_2\geq k'$.
\kcl
\dow
%Suppose the contrary. The number of uncrossed squares is equal to $w'=k'd'$.
The columns of index bigger than $o'$ have at most $y_2$ uncrossed boxes, and there
are $l'$ of them. So the total number of uncrossed boxes is bounded by $y_1+l'y_2$.
This number is also $k'd'$, and since $l'\le d'$ we get
 $k'd'\leq y_1+l'y_2\leq d'(y_1+y_2)$, hence $k'\le y_1+y_2$.
\kdow

\medskip
Now suppose that step C produces a repetition when the two permutations
involved are the identity. This would mean that we cross two boxes in $\lambda$
belonging to the same row. But then, consider the situation at the step just before.
Since we are about to apply step C we have $l'<d'$.
Define $y_1$ and $y_2$ as above. Then $y_2\le\lambda_{o'}-m'<k'$, since each
uncrossed entry in the column $o'+1$ is either going to be crossed or is at the
same height of a box in column $o'$ that is going to be crossed, and by hypothesis
there is some row at which both events will occur. But then the total number
of uncrossed boxes is bounded by $(l'+1)(k'-1)<d'k'$, a  contradiction!

\medskip\noindent {\bf The contribution of $Q$ is positive}.
We only have to prove that if a simple tensor in $P$ gives $Q$ after symmetrization,
it has to come with a positive sign. First observe that we can always suppose
that the permutations giving such a simple tensor are increasing on the
set of indices contributing to the same column of our rectangular tableaux.
Otherwise we can rearrange them and get the same contribution. The main
observation is that once this hypothesis has been made, these permutations
are necessarily {\it paired}.

\dfi[Paired set, paired permutation]
We say that a set of integers is paired if whenever it contains $i$ odd,
it also contains $i+1$.
We say that a permutation $\sigma$ is paired if $\sigma(j)=i$ odd
implies $\sigma(j+1)=i+1$.
\kdfi
This implies that each such $j$ is odd. Therefore any paired permutations
has positive sign.
%Moreover sum and difference of
%paired sets is a paired set.
We observe that for each wedge product $e_{a_1}\wedge\dots\wedge
e_{a_k}$ appearing in $Q$ the set $\{a_1,\dots,a_k\}$ is paired.
Indeed, since $k$ is even, in each wedge product defining $P$
each permutation appears an even number of times. Hence when these
permutations are all equal to the identity, the indices in these wedge
products form paired sets. Also the frozen indices appearing in
each wedge product form paired sets.

Suppose that a simple tensor $S$ in $P$ has symmetrization $Q$. Then the indices
appearing in each wedge product in $S$ must form a paired set.
We deduce inductively that the indices in each permutation
in each wedge product form a paired set.
Therefore each of these permutations  is itself paired. But then their signs
are positive, hence they must contribute positively to $Q$. \qed


\begin{thebibliography}{aa}

\bibitem[Br]{br}
Brion M., {\it Plethysm and Verma modules}
J. London Math. Soc. {\bf 52} (1995), 449--466.

%Brion M., {\it Stable properties of plethysm: on two conjectures of Foulkes},
%Manuscripta Math. {\bf 80} (1993), 347–371.

\bibitem[BCI]{bci}
B\"urgisser P., Christandl M.,  Ikenmeyer C., {\it
Even partitions in plethysms}, J. Algebra {\bf 328} (2011), 322--329.

\bibitem[FH]{fh}
Fulton W., Harris J., Representation Theory - A First Course,
Grad. Texts in Math. {\bf 129}, Springer 1991.

\bibitem[LaR]{lar} Langley T., Remmel J., {\it
The plethysm $s_\lambda [s_\mu ]$ at hook and near-hook shapes},
Electron. J. Combin. {\bf 11} (2004), no. 1, Research Paper 11, 26 pp.

\bibitem[LoR]{lor} Loehr N., Remmel J., {\it
A computational and combinatorial expos\'e of plethystic calculus},
J. Algebraic Combin. {\bf 33} (2011), 163--198.

\bibitem[Mc]{mcd}
Macdonald I.,
Symmetric functions and Hall polynomials, Oxford Mathematical Monographs,
Second edition, Oxford University Press 1995.

\bibitem[Ma]{man}
Manivel L., {\it Gaussian maps and plethysm}, in Algebraic geometry
(Catania, 1993/Barcelona, 1994), 91--117,
Lecture Notes in Pure and Appl. Math. {\bf 200}, Dekker  1998.

\bibitem[We]{we}
Weintraub  S., {\it Some observations on plethysms},
J. Algebra {\bf 129} (1990), 103--114.

\end{thebibliography}
\end{document}